\documentclass{birkmult}
\begin{document}
%
%
%
 \newtheorem{thm}{Theorem}[section]
 \newtheorem{cor}[thm]{Corollary}
 \newtheorem{lem}[thm]{Lemma}
 \newtheorem{prop}[thm]{Proposition}
 \theoremstyle{definition}
 \newtheorem{defn}[thm]{Definition}
\theoremstyle{assumption}
 \newtheorem{assumption}[thm]{Assumption}
 \theoremstyle{remark}
 \newtheorem{rem}[thm]{Remark}
 \newtheorem*{ex}{Example}
 \numberwithin{equation}{section}
 \def\R{\mathbb{R}}
\def\N{\mathbb{N}}
\def\K{\mathbf{K}}
\def\Q{\mathbb{Q}}
\def\C{\mathbb{C}}
\def\P{\mathbb{P}}
\def\Z{\mathbb{Z}}
\def\S{\mathbb{S}}

%
%
%
%
%
%
%
%
%
\title{Representation of nonnegative convex polynomials}

\author[Lasserre]{Jean B. Lasserre}

\address{%
LAAS-CNRS and Institute of Mathematics\\
University of Toulouse\\
LAAS, 7 avenue du Colonel Roche\\
31077 Toulouse C\'edex 4\\
France}

\email{lasserre@laas.fr}

\thanks{This work was completed with the support of 
the (french) ANR grant NT05-3-41612.}



\date{}

\begin{abstract}
We provide a specific representation of
convex polynomials nonnegative on a convex (not necessarily compact) 
basic closed semi-algebraic set $\K\subset\R^n$. Namely, 
they belong to a specific subset of the quadratic module generated by the concave polynomials that define $\K$. 
\end{abstract}

\maketitle

\section{Introduction}
An important research area of real algebraic geometry is concerned with representations of polynomials positive on a basic semi-algebraic set 
\begin{equation}
\label{setk}
\K:=\{x\in\R^n\::\:g_j(x)\geq0,\quad j=1,\ldots,m\}\subset\R^n
\end{equation}
where $g_j\in\R[X]$, $j=1,\ldots,m$.

An important result in this vein is Schm\"udgen's Positivstellensatz \cite{schmudgen} which states that if $\K$ is compact and
$f\in\R[X]$ is positive on $\K$ then $f$ belongs to the preordering $P(g)$
generated by the $g_j$'s; bounds on the degrees
in the representation are even provided in Schweighofer \cite{markus}.
Under a rather weak additional assumption on the $g_j$'s,
Putinar's refinement \cite{putinar} states that
$f$ even belongs to the quadratic module $Q(g)$ generated by the $g_j$'s.
The  above mentioned representation
results do not specialize when $f$ is convex
and the $g_j$'s are concave (so that $\K$ is convex)
a highly important case, particularly in optimization.
Also, as soon as $\K$ is not compact any more then negative results, notably by Scheiderer \cite{claus}, exclude
to represent {\it any} $f$ positive on $\K$ as an element of
$P(g)$ or $Q(g)$ (except perhaps in low-dimensional cases). For more details,
the interested reader is referred to the nice survey \cite{claus}.

However, inspired and motivated by some classical results from convex optimization, we show that specialized representation results are possible
when $f$ is convex and the $g_j$'s are concave, in which case
$\K\subset\R^n$ is a closed (not necessarily compact) convex basic 
semi-algebraic set. Namely, a specific subset $Q_c(g)$ of the quadratic module $Q(g)$
is such that $Q_c(g)\cap F$ is {\it dense} (for the $l_1$-norm of coefficients) in the convex cone $F$ of convex polynomials, nonnegative on $\K$.

\section{Convex polynomials on a convex semi-algebraic set}
\subsection{Notation and Preliminaries}
Let $\R[X]$ be the ring of real polynomials in the variables $X=(X_1,\ldots,X_n)$, and let $\Sigma^2\subset\R[X]$ be the subset of sums of squares (sos) polynomials. If $f\in\R[X]$, write
$f(X)=\sum_{\alpha\in\N^n}f_\alpha X^\alpha$, and denote
its $l_1$-norm by $\Vert f\Vert_1 \,(=\sum_{\alpha\in\N^n}\vert f_\alpha\vert$).

Let $Q(g)\subset\R[X]$ 
be the {\it quadratic module} generated by a set
of polynomials $g=(g_j)_{j=1}^m\subset\R[X]$, that is,
\begin{equation}
\label{quad}
Q(g)\,:=\,\left\{\: \sigma_0+\sum_{j=1}^m\sigma_j\,g_j\::\quad \sigma_j\in\Sigma^2,\:j=0,\ldots,m\:\right\}.
\end{equation}
Throughout the paper we make the following assumption.
\begin{assumption}
\label{ass1}
$\K\subset\R^n$ is defined in (\ref{setk}) and is such that:

{\rm (a)} $g_j$ is concave for every $j=1,\ldots,m$.

{\rm (b)} There exists $z\in\K$ such that
$g_j(z)>0$ for every $j=1,\ldots,m$.
\end{assumption}
Assumption \ref{ass1}(b), known as Slater condition,
is an important regularity condition for the celebrated Karush-Kuhn-Tucker optimality conditions.
\begin{prop}
\label{prop1}
Let Assumption \ref{ass1} hold
and let $f\in\R[X]$ be convex and such that
$f^*:=\inf_x\{f(x)\::\:x\in\K\}=f(x^*)$ for some $x^*\in\K$.

Then there exists $\lambda\in\R^m_+$ such that
\begin{equation}
\label{kkt}
\nabla f(x^*)-\sum_{j=1}^m\lambda_j\nabla g_j(x^*)\,=\,0\,;\quad
\lambda_j\,g_j(x^*)\,=\,0,\quad j=1,\ldots,m.
\end{equation}
In other words, the Lagrangian $L_f\in\R[X]$ defined  by
\begin{equation}
\label{lagrangian}
X\mapsto L_f(X)\,:=\,f(X)-f^*-\sum_{j=1}^m\lambda_j\,g_j(X),\qquad X\in\R^n,\end{equation}
is a nonnegative polynomial which satisfies
\begin{equation}
\label{kkt2}
L_f(x^*)\,=\,0\,;\quad \nabla L_f(x^*)=0.
\end{equation}
\end{prop}
See e.g. Polyak \cite{polyak}.

\subsection{Convex Positivstellensatz}

If one is interested in representation of polynomials nonnegative on $\K$,
the first polynomial to consider is of course $f-f^*$ where
$0\leq f^*=\inf_{x\in\K}f(x)$. Indeed, any other positive polynomial
is just $(f-f^*)+f^*$ with $f^*\geq0$. 
And so, if $f-f^*$ belongs to some preordering or some quadratic module, then so does $f$.
From Proposition \ref{prop1} it is easy to establish the following result.
\begin{cor}
\label{coro1}
Let Assumption \ref{ass1} hold and 
let $f\in\R[X]$ be convex and such that
$f^*:=\inf_x\{f(x)\::\:x\in\K\}=f(x^*)$ for some $x^*\in\K$.
If the nonnegative polynomial 
$L_f$ of (\ref{lagrangian}) is sos then 
\begin{equation}
\label{main}
f-f^*\,=\,\sigma +\sum_{j=1}^n\lambda_j g_j
\end{equation}
for some convex sos polynomial
$\sigma\in\Sigma^2$ and some nonnegative scalars $\lambda_j$, $j=1,\ldots,m$. That is, 
$f-f^*\in Q(g)$, with $Q(g)$ as in (\ref{quad}).
In addition, 
the sos weights associated with the $g_j$'s are just nonnegative constants, and $\sigma$ is convex.
\end{cor}
\begin{proof}
Follows from the definition (\ref{lagrangian}) of $L_f$, and the fact that $L_f$ is sos.
\end{proof}
Hence in view of Corollary \ref{coro1}, an interesting issue is to provide sufficient conditions for $L_f$ to be sos. For instance,
consider the following definition from 
Helton and Nie \cite{helton}
\begin{defn}[Helton and Nie \cite{helton}]
\label{def1}
A polynomial $f\in\R[X]$ is sos-convex if its Hessian $\nabla ^2f$ is a sum of squares (sos), that is, there is some integer $p$ and
some matrix polynomial $F\in\R[X]^{p\times n}$
such that
\begin{equation}
\label{sos}
\nabla^2f(X)\,:=\,\left(\frac{\partial^2f(X)}{\partial X_i\partial X_j}\right)_{ij}\,=\,F(X)^TF(X).
\end{equation}
\end{defn}
\begin{cor}
\label{coro2}
Let Assumption \ref{ass1} hold, and let $f\in\R[X]$ be convex and such that
$f^*:=\inf_x\{f(x)\::\:x\in\K\}=f(x^*)$ for some $x^*\in\K$.

If $f$ is sos-convex and $-g_j$ is sos-convex
for every $j=1,\ldots,m$, then $f-f^*\in Q(g)$. More precisely,
(\ref{main}) holds for some convex sos polynomial
$\sigma\in\Sigma^2$ and some nonnegative scalars $\lambda_j$, $j=1,\ldots,m$. 
\end{cor}
\begin{proof}
From Proposition \ref{prop1}, let $L_f$ be as in (\ref{lagrangian}).
As $f$ and $-g_j$ are sos convex, write
\[\nabla^2f(X)\,=\,F(X)^TF(X);\quad
-\nabla^2g_j(X)\,=\,G_j(X)^TG_j(X),\quad j=1,\ldots,m,\]
for some $F\in\R[X]^{p\times n}$ and some
$G_j\in\R[X]^{p_j\times n}$, $j=1,\ldots,m$. Hence,
\[\nabla^2L_f\,=\,\nabla^2f-\sum_{j=1}^m\lambda_j\nabla^2g_j\,=\,
F^TF+\sum_{j=1}^m\lambda_jG_j^TG_j\,=\,H^TH,\]
with $H^T:=\left[F^T\,\vert\,\sqrt{\lambda_1}\,G_1^T\vert\cdots\vert\,\sqrt{\lambda_m}\,G_m^T\right]$,
and so $L_f$ is sos-convex.
As (\ref{kkt2}) holds, by Lemma 3.2 in Helton and Nie \cite{helton}, the polynomial $L_f$ is sos, and so, by Corollary \ref{coro1}, the desired result (\ref{main}) holds.
\end{proof}
Next, consider the subset $Q_c(g)\subset Q(g)$ defined by:
\begin{equation}
\label{qc}
Q_c(g)\,:=\,\left\{\:\sigma+\sum_{j=1}^m\lambda_j\,g_j\::\quad
\lambda\in\R^m_+\,;\:\sigma\in\Sigma^2,\:\sigma\mbox{ convex.}
\:\right\}\subset Q(g).
\end{equation}
The set $Q_c(g)$ is a specialization of
$Q(g)$ to the convex case, in that the weights asociated with the $g_j$'s are nonnegative scalars, i.e., sos polynomials of degree 0, and the sos polynomial $\sigma$ is convex.
\begin{thm}
\label{thmain2}
Let Assumption \ref{ass1} hold, and let $Q_c(g)$ be as in (\ref{qc}).
Let $F\subset\R[X]$ be the convex cone
of convex polynomials nonnegative on $\K$.

Then $Q_c(g)\cap F$ is dense in $F$ for the $l_1$-norm $\Vert \cdot\Vert_1$. In particular, 
if $\K=\R^n$ (so that $F$ is now the set of nonnegative convex polynomials), then
$\Sigma^2\cap F$ is dense in $F$.
\end{thm}
\begin{proof}
Let $f\in F$ and let $r_0:=\lfloor ({\rm deg}\,f)/2\rfloor+1$.
Given $r\in\N$, let $\Theta_r\in\R[X]$ be the polynomial
\begin{equation}
\label{Theta}
X\mapsto \Theta_r(X)\,:=\,1+\sum_{i=1}^nX_i^{2r}.
\end{equation}
For every $\epsilon>0$, the polynomial
$f_{\epsilon 0}(X):=f(X)+\epsilon\,\Theta_{r_0}(X)$ is convex
and nonnegative on $\K$, i.e., $f_{\epsilon 0}\in F$. In addition, 
\[0\,\leq\,f^*:=\inf_{x\in\K}f(x)\,\leq\,\inf_{x\in\K}f_{\epsilon 0}(x)
\,=\,f_{\epsilon 0}(x^*_\epsilon)\,=:\,
f_{\epsilon}^*,\]
for some $x^*_\epsilon\in\K$. 
Indeed, the level set
$\{x\in\K\::\:f_{\epsilon 0}(x)\leq \alpha\}$ 
is compact for every $\alpha\in\R$, and so,
$f_{\epsilon 0}$ attains its minimum on $\K$.
Obviously, we also have
$\Vert f_{\epsilon 0}-f\Vert_1\to0$
as $\epsilon\downarrow 0$.
Next, let $L_{f_{\epsilon 0}}$ 
be as in (\ref{lagrangian}), i.e.,
\[L_{f_{\epsilon 0}}\,=\,
f+\epsilon\,\Theta_{r_0}-f^*_\epsilon-\sum_{j=1}^n
\lambda_j^\epsilon\,g_j,\]
for some nonnegative vector $\lambda_j^\epsilon\in\R^m_+$.
As $L_{f_{\epsilon 0}}\geq0$ on $\R^n$,
by Corollary 3.3 in Lasserre and Netzer \cite{lassnetzer},
there exists $r_\epsilon\in\N$ such that for every
$r\geq r_\epsilon$, $L_{f_{\epsilon 0}}+\epsilon\,\Theta_r$ is sos.
That is, $\sigma:=L_{f_{\epsilon 0}}+\epsilon\,\Theta_r\in\Sigma^2$ and so
\[f_\epsilon\,:=\,f+\epsilon\,(\Theta_{r_0}+\Theta_r)\,=\,\sigma+f^*_\epsilon+\sum_{j=1}^n\lambda_j^\epsilon\,g_j.\]
Notice that by definition, $\sigma\in\Sigma^2$ is convex.
Next, as $f^*_\epsilon\geq0$, $\sigma+f^*_\epsilon\in\Sigma^2$, 
and so, equivalently, $f_\epsilon\in Q_c(g)$. 

In addition, $f_\epsilon\in F$ because
$f_\epsilon$ is convex (as $f_\epsilon=
f+\epsilon (\Theta_{r_0}+\Theta_r)$) and nonnegative on $\K$
(as $f_\epsilon\geq f$), and so,
$f_\epsilon\in Q_c(g)\cap F$.
Finally, $\Vert f-f_\epsilon\Vert_1=\epsilon\Vert\Theta_{r_0}+\Theta_r\Vert_1\to 0$ as $\epsilon\downarrow 0$.

Finally, if $\K=\R^n$ (so that $F$ is now the set of nonnegative convex polynomials),
one obtains $Q_c(g)=\Sigma^2$.
\end{proof}
One may also replace $\Theta_r$ in (\ref{Theta}) with the new perturbation
\[X\mapsto \theta_r(X)\,:=\,\sum_{k=0}^r\sum_{j=1}^n\frac{X_j^{2k}}{k{\rm !}}.\]
This perturbation also preserves convexity.
In addition, not only $\Vert f-f_\epsilon\Vert_1\to 0$
as $\epsilon\downarrow 0$, but 
the convergence $f_\epsilon\to f$ is also
{\it uniform} on compact sets!

\end{document}